\documentclass[10pt]{article}
\usepackage{chicago,epic,graphpap,graphics,float,tabularx,paralist,longtable,fancyhdr,fancyvrb}
\usepackage{amsmath}
\usepackage{amstext}
\usepackage{amsbsy}
\usepackage{amsthm}
\usepackage{amsfonts}
\usepackage{amssymb}
\usepackage{amscd}
\usepackage{eucal}
\usepackage{bbm}
\def\X{{\widehat X}}
\def\XX{{\widetilde X}}
\def\Y{{\widehat Y}}

\def\ito{It{\^o}}
\def\as{\quad\text{{\rm a.s.}}}
\def\R{{\mathbb R}}
\def\eqdef{\triangleq}
\def\half{\frac{1}{2}}
\def\sumi{\sum_{i=1}^n}
\def\sumj{\sum_{j=1}^n}
\def\sumk{\sum_{k=1}^n}
\def\brac#1{\langle #1 \rangle}

\def\dd{\circ d}
\def\F{{\mathcal F}}

\def\T{{\mathcal T}}
\def\eps{\varepsilon}

\def\intt{\int_0^t}

\def\1{{\mathbbm 1}}

\def\s{\sigma}
\def\sgn{{\rm sgn}}
\def\m{\mu}
\def\p{\pi}
\def\D{\Delta}
\def\T{{\cal T}}
\def\S{{\bf S}}

\begin{document}

\centerline{\Large\bf  Stratonovich representation of semimartingale rank processes}
\vspace{10pt} \centerline{\large  Robert Fernholz\footnote{INTECH, One Palmer Square, Princeton, NJ 08542.  bob@bobfernholz.com. The author thanks Ioannis Karatzas and Mykhaylo Shkolnikov for their invaluable comments and suggestions regarding this research.}} \centerline{\today}

\begin{abstract}
Suppose that $X_1,\ldots,X_n$ are continuous semimartingales that are reversible and have nondegenerate crossings. Then the corresponding rank processes can be represented by generalized Stratonovich integrals, and this  representation can be used to decompose the relative log-return of portfolios generated by functions of ranked market weights.
\end{abstract}

\vspace{10pt}
\noindent{\large\bf Introduction}
\vspace{10pt}

For $n\ge2$, consider a family of continuous semimartingales $X_1,\ldots,X_n$  defined on $[0,T]$ under the usual filtration $\F^X_t$, with quadratic variation processes $\brac{X_i}$.  Let $r_t(i)$ be the rank of $X_i(t)$, with $r_t(i)<r_t(j)$ if $X_i(t)>X_j(t)$ or if $X_i(t)=X_j(t)$ and $i<j$. The corresponding {\em rank processes} $X_{(1)},\ldots,X_{(n)}$  are defined by $X_{(r_t(i))}(t)=X_i(t)$. We shall show that if the $X_i$ are {\em reversible} and have {\em nondegenerate crossings,}  then the rank processes can be represented by
\begin{equation}\label{1}
dX_{(k)}(t)=\sumi\1_{\{X_i(t)=X_{(k)}(t)\}}\dd  X_i(t),\as,
\end{equation}
where $\hspace{-7pt}\phantom{a}\dd$ is the generalized Stratonovich integral developed by \citeN{Russo:2007}.

An  {\em Atlas model} is a family of positive continuous semimartingales $X_1,\ldots,X_n$ defined as an \ito\ integral on $[0,T]$ by 
\begin{equation}\label{0}
d\log X_i(t)=\big(-g+ng\1_{\{r_t(i)=n\}}\big)dt+\s\,dW_i(t),
\end{equation}
 where $g$ and $\s$ are positive constants and $(W_1,\ldots,W_n)$ is a Brownian motion (see \citeN{F:2002}). Here the $X_i$ represent the capitalizations of the companies in a stock market, and $d\log X_i$ represents the log-return of the $i$th stock.  We shall show the representation~\eqref{1} is valid for the Atlas rank processes $\log X_{(k)}$.

In \citeN{F:2016} it was shown that in a stock market with stocks represented by positive continuous semimartingales, under certain conditions the log-return of a portfolio can be decomposed into a {\em structural process} and a {\em trading process,} and for a portfolio generated by a $C^2$ function of the market weight processes, these components correspond to the log-change in the generating function and the drift process (see \citeN{F:rank}). The Stratonovich representation~\eqref{1} allows us to extend this decomposition to portfolios generated by $C^2$ functions of the ranked market weight processes in Atlas models. 
 
\vspace{15pt}
\noindent{\large\bf \ito\ integrals and Stratonovich integrals}
\vspace{10pt}
 
 Let $X$ and $Y$ be continuous semimartingales on $[0,T]$ with the filtration $\F_t^{X,Y}$. Then the {\em Fisk-Stratonovich integral} is defined by
 \begin{equation}\label{3}
 \intt Y(s)\dd X(s) \eqdef \intt Y(s)\,dX(s) + \half\brac{Y,X}_t,
 \end{equation}
for $t\in[0,T]$, where the integral on the right hand side is the \ito\ integral and $\brac{X,Y}_t$ is the cross variation of $X$ and $Y$ over $[0,t]$ (see \citeN{Karatzas/Shreve:1991}). The Fisk-Stratonovich integral is defined only for  semimartingales, but in some cases can be extended to more general integrands.  Following \citeN{Russo:2007},  Definition~1, for a continuous semimartingale $X$  and a locally integrable process $Y$, both defined on $[0,T]$, we define the {\em forward integral, backward integral,} and {\em covariation process} by
{\allowdisplaybreaks
\begin{align}
\intt Y(s)\,d^-X(s)&\eqdef \lim_{\eps\downarrow0}\intt Y(s)\frac{X(s+\eps)-X(s)}{\eps}ds\label{0.1}\\
\intt Y(s)\,d^+X(s)&\eqdef \lim_{\eps\downarrow0}\intt Y(s)\frac{X(s)-X(s-\eps)}{\eps}ds\label{0.2}\\
[X,Y]_t&\eqdef \lim_{\eps\downarrow0}\intt \frac{(X(s+\eps)-X(s))(Y(s+\eps)-Y(s))}{\eps}ds,\label{0.3}
\end{align}
}
for $t\in[0,T]$, where the limits are uniform in probability on $[0,T]$.  We shall use the convention of \citeN{Russo:2007} that for the evaluation of these limits a continuous function $X$ defined on $[0,T]$ is implicitly extended to $\R$ by setting $X(t)=X(0)$ for $t<0$ and $X(t)=X(T)$ for $t>T$. Then, by \citeN{Russo:2007}, Definition~10, the {\em Stratonovich integral} is given by
  \begin{equation}\label{3.1}
 \intt Y(s)\dd X(s) \eqdef \intt Y(s)\,dX(s) + \half\big[Y,X\big]_t,
 \end{equation}
where the integral on the right hand side is the \ito\ integral. If both $X$ and $Y$ are continuous semimartingales, then
\[
\brac{X,Y}_t=\big[X,Y\big]_t,\as,
\]
and the Stratonovich integral is equivalent to the Fisk-Stratonovich integral.

For a continuous semimartingale $X$ and $C^2$ function $F$ defined on the range of $X$, \ito's rule establishes that
\[
F(X(t))-F(X(0))= \intt F'(X(s))\,dX(s)+\half \intt F''(X(s))\,d\brac{X}_s,\as,
\]
and with the Fisk-Stratonovich integral, this becomes 
\begin{equation}\label{4}
F(X(t))-F(X(0))=\intt F'(X(s))\dd X(s),\as,
\end{equation}
as in ordinary calculus (see \citeN{Karatzas/Shreve:1991}). The relationship \eqref{4} can be extended to a wider class of functions in some cases. For example, for an absolutely continuous function $F$ and Brownian motion $W$, it was  shown in \shortciteN{Follmer:1995}, Corollary~4.2, that \eqref{4} holds, 
so for the absolute-value function we have
\begin{equation}\label{20}
|W(t)|=\intt\sgn(W(s))\dd W(s),\as,
\end{equation}
where $\sgn(x)\eqdef \1_{\{x>0\}}-\1_{\{x\le0\}}$. The \citeN{Russo:2007} results allow us to extend this relationship to a class of continuous semimartingales.

\vspace{15pt}
\noindent{\bf Definition~1.}  Let $X$ be a continuous semimartingale defined on $[0,T]$ under the filtration $\F^X_t$. Then $X$ is {\em   reversible} if  the time-reversed process $\X$ defined by $\X(t)=X(T-t)$ is also a continuous semimartingale on $[0,T]$ under the time-reversed filtration $\F^\X_t$. 

\vspace{15pt}
\noindent{\bf Definition~2.} The continuous semimartingales $X_1,\ldots,X_n$ have {\em   nondegenerate crossings}  if  for any $i\ne j$ the set $\{t:X_i(t)= X_j(t)\}$ almost surely has measure zero with respect to $d\brac{X_k}_t$, for $1\le k\le n$.

\vspace{15pt}
\noindent{\bf Lemma~1.}  {\em Suppose that the continuous semimartingales $X_1,\ldots,X_n$ have nondegenerate crossings. Then the same is true for the rank processes $X_{(1)},\ldots,X_{(n)}$.}

\vspace{10pt}
\noindent{\em Proof.} Let $p_t\in\Sigma_n$ be the inverse permutation to the rank function $r_t$. If $X_{(k)}(t)=X_{(\ell)}(t)$ for $k\ne\ell$, then $X_i(t)=X_j(t)$ for $i=p_t(k)\ne p_t(\ell)=j$, so 
\begin{equation}\label{15}
\bigcup_{i\ne j}\{t:X_i(t)=X_j(t)\}=\bigcup_{k\ne\ell}\{t:X_{(k)}(t)=X_{(\ell)}(t)\}. 
\end{equation}

From \citeN{BG:2008} we have the representation
\begin{equation*}
dX_{(k)}(t)=\sumi\1_{\{X_i(t)=X_{(k)}(t)\}}dX_i(t)+\text{ finite variation terms},\as,
\end{equation*}
for $k=1,\ldots,n$, so $d\brac{X_{(k)}}_t\ll d\brac{X_1}_t+\cdots+d\brac{X_1}_t$, for $k=1,\ldots,n$. In the same manner we can show that  $d\brac{X_i}_t\ll d\brac{X_{(1)}}_t+\cdots+d\brac{X_{(n)}}_t$, for $i=1,\ldots,n$, so sets of the form $\{t:X_i(t)= X_j(t)\}$, for $i\ne j,$ or $\{t:X_{(k)}(t)= X_{(\ell)}(t)\}$, for $k\ne\ell$, will almost surely have measure zero with respect to $d\brac{X_i}_t$, for $i=1,\ldots,n$, and with respect to $d\brac{X_{(k)}}_t$, for $k=1,\ldots,n$, and the same holds for finite unions of such sets, as in \eqref{15}.\qed

\vspace{15pt}
\noindent{\bf Lemma~2.}  {\em Let $X$ be a reversible continuous semimartingale defined on $[0,T]$, and suppose that the set $\{t:X(t)=0\}$ almost surely has measure zero with respect to $d\brac{X}_t$.  Then
\begin{equation}\label{16}
|X(t)|-|X(0)|=\intt\sgn(X(s))\dd X(s),\as
\end{equation}
}

\noindent{\em Proof.} References in this proof denoted by R\&V are from  \citeN{Russo:2007}.

The Tanaka-Meyer formula states that for the \ito\ integral
\begin{equation}\label{16.1}
\intt\sgn(X(s))\,d X(s)= |X(t)|-|X(0)|+2\Lambda_X(t),\as,
\end{equation}
where  $\Lambda_X$ is the local time at zero for $X$ (see \citeN{Karatzas/Shreve:1991}). From R\&V, Proposition~1, 
\begin{equation}\label{16.2}
\intt\sgn(X(s))\dd X(s)= \half\bigg(\intt\sgn(X(s))\,d^- X(s)+\intt\sgn(X(s))\,d^+ X(s)\bigg),
\end{equation}
with the forward and backward integrals defined by \eqref{0.1} and \eqref{0.2}. Since $\sgn(X)$ is continuous outside the set $\{t:X(t)=0\}$, and this set almost surely has measure zero with respect to $d\brac{X}_t$, R\&V Proposition~6 implies that
\begin{align}
\intt\sgn(X(s))\,d^- X(s)&=\intt\sgn(X(s))\,d X(s)\notag\\
&=|X(t)|-|X(0)|+2\Lambda_X(t),\as,\label{16.3}
\end{align}
by equation \eqref{16.1}.

By R\&V, Proposition~1, 
\begin{equation}\label{16.4}
\intt\sgn(X(s))\,d^+X(s)=-\int_{T-t}^T\sgn(\X(s))\,d^- \X(s),
\end{equation}
where $\X$ is the time-reversed version of $X$. By hypothesis, $\X$ is a continuous semimartingale  on $[0,T]$ with respect to the reverse filtration, so as in \eqref{16.3} we have
\begin{align}
\int_{T-t}^T\sgn(\X(s))\,d^- \X(s)&=|\X(T)|-|\X(T-t)|+2\big(\Lambda_\X(T)-\Lambda_\X(T-t)\big)\notag\\
&=|X(0)|-|X(t)|+2\Lambda_X(t),\as\label{16.5}
\end{align}
If we combine \eqref{16.2}, \eqref{16.3}, \eqref{16.4}, and \eqref{16.5}, then \eqref{16} follows.\qed

\vspace{15pt}
\noindent{\bf Lemma~3.}  {\em Let $ X$ and $ Y$ be reversible continuous semimartingales defined on $[0,T]$ under under a common filtration and suppose that they have nondegenerate crossings. Then
\begin{multline}\label{6.5}
\big| X(t)- Y(t)\big|-\big| X(0)- Y(0)\big|=\intt\sgn\big( X(s)-Y(s)\big)\dd  X(s)\\ -\intt\sgn\big( X(s)-Y(s)\big)\dd Y(s),\as
\end{multline}
}

\noindent{\em Proof.} References in this proof denoted by R\&V are from  \citeN{Russo:2007}.

Since $ X-Y$ is a   reversible continuous semimartingale, it follows from Lemma~2 that
\[
\big| X(t)-Y(t)\big|-\big| X(0)-Y(0)\big|=\intt\sgn\big( X(s)-Y(s)\big)\dd \big( X(s)-Y(s)\big),\as,
\]
so, due to linearity of the integral with respect to the differentials, it suffices to show that the integrals in \eqref{6.5} are defined. Let us first consider the integral with respect to $d X$.

By the definition of the Stratonovich integral in \eqref{3.1}, 
\begin{equation}\label{6.7}
\intt\sgn\big( X(s)-Y(s)\big)\dd  X(s) = \intt\sgn\big( X(s)-Y(s)\big)d  X(s) +\big[\sgn\big( X-Y\big),  X\big]_t,
\end{equation}
if the terms on the right hand side are defined. The \ito\ integral in \eqref{6.7} is defined, so we need only  consider the covariation term.  From R\&V Proposition~1, 
\begin{equation}\label{6.6}
\big[\sgn\big( X-Y\big),  X\big]_t=\intt \sgn\big( X(s)-Y(s)\big)d^+  X(s)-\intt\sgn\big( X(s)-Y(s)\big)d^-  X(s),
\end{equation}
and will be defined if the two integrals are. Since $\sgn( X-Y)$ is continuous outside $\{t: X(t)=Y(t)\}$, which almost surely has measure zero with respect to $d\brac{ X}_t$,  R\&V  Proposition~6 implies that 
\begin{equation}\label{6.9}
\intt\sgn\big( X(s)-Y(s)\big)d^-  X(s) = \intt\sgn\big( X(s)-Y(s)\big)d  X(s),
\end{equation}
and since this \ito\ integral is defined, so is the forward integral. By R\&V  Proposition~1, 
\begin{equation}\label{6.8}
\intt\sgn\big( X(s)-Y(s)\big)d^+  X(s)=-\int_{T-t}^T\sgn\big(\X(s)-\Y(s)\big)d^- \X(s),
\end{equation}
where $\X$ and $\Y$ are the time-reversed versions of $X$ and $Y$ on $[0,T]$. By hypothesis, the time-reversed process $\X$ is a continuous semimartingale, so as in \eqref{6.9} its forward integral is defined, and this defines the backward integral in \eqref{6.8}. Hence, the covariation in \eqref{6.6} is defined, so both terms on the right hand side of \eqref{6.7} are defined, and this defines the Stratonovich integral with respect to $d X$ in \eqref{6.5}.

The same reasoning holds for the integral with respect to $dY$. \qed

\vspace{15pt}
\noindent{\large\bf A Stratonovich representation for rank processes}
 \vspace{10pt}

We would like to prove \eqref{1}, and we shall start with a lemma that establishes this result for $n=2$ and then apply the lemma to prove the general case with $n\ge2$.

\vspace{10pt}
\noindent{\bf Lemma~4.}  {\em Let $X_1$ and $X_2$ be reversible continuous semimartingales defined on $[0,T]$ under a common filtration,  and suppose that they have nondegenerate crossings. Then
\begin{equation}\label{8.0}
X_1(t)\lor X_2(t)-X_1(0)\lor X_2(0)= \intt\1_{\{X_1(s)\ge X_2(s)\}}\dd X_1(s)+\intt\1_{\{X_1(s)< X_2(s)\}}\dd X_2(s),\as,
\end{equation}
and
\begin{equation}\label{8.1}
X_1(t)\land X_2(t)-X_1(0)\land X_2(0)= \intt\1_{\{X_1(s)< X_2(s)\}}\dd X_1(s)+\intt\1_{\{X_1(s)\ge X_2(s)\}}\dd X_2(s),\as
\end{equation}
}

\noindent{\em Proof.} For $t\in[0,T]$ we have
\[
  X_1(t)\lor   X_2(t)= \half\Big(  X_1(t)+  X_2(t)+\big|  X_2(t)-  X_1(t)\big|\Big),\as,
\]
so by Lemma~3,
{\allowdisplaybreaks
\begin{align*}
X_1(t)\lor   X_2(t)-X_1(0)\lor   X_2(0)
&= \half\Big(X_1(t)+ X_2(t)-X_1(0)- X_2(0)+\intt\sgn\big(X_2(s)-X_1(s)\big)\dd   X_2(s)\\
&\hspace{80pt}-\intt\sgn\big( X_2(s)-X_1(s)\big)\dd   X_1(s)\Big)\\
&= \half\intt \Big(1+\sgn\big(X_2(s)-X_1(s)\big)\Big)\dd   X_2(s)\\
&\hspace{80pt}+\half\intt \Big(1-\sgn\big( X_2(s)-X_1(s)\big)\Big)\dd   X_1(s)\\
&=\intt \1_{\{X_1(s)<X_2(s)\}}\dd  X_2(s)+\intt\1_{\{X_1(s)\ge X_2(s)\}}\dd   X_1(s),\as,
\end{align*}}
which proves \eqref{8.0}. Equation \eqref{8.1} follows from this and the fact that
\[
X_1(t)\land   X_2(t)= X_1(t)+X_2(t)-  X_1(t)\lor   X_2(t),\as\tag*{\qed}
\]

\noindent{\bf Proposition~1.} {\em  Let  $X_1,\ldots,X_n$ be continuous semimartingales  defined on $[0,T]$ that are reversible and have nondegenerate crossings.  Then the rank processes $X_{(1)},\ldots,X_{(n)}$ satisfy}
\begin{equation}\label{11}
dX_{(k)}(t)=\sumi\1_{\{X_i(t)=X_{(k)}(t)\}}\dd  X_i(t),\as
\end{equation}

\vspace{10pt}

\noindent{\em Proof.} It follows from Lemma~4 that \eqref{11} holds for $n=2$, so let us assume that it holds for $X_1,\ldots,X_{n-1}$, and prove that it then holds for $X_1,\ldots,X_n$. Let $\XX_{(1)},\ldots,\XX_{(n-1)}$ be the ranked processes $X_1,\ldots,X_{n-1}$, so by our inductive hypothesis we have
\begin{equation}\label{9}
d\XX_{(k)}(t)=\sum_{i=1}^{n-1}\1_{\{X_i(t)=\XX_{(k)}(t)\}}\dd  X_i(t),\as,
\end{equation}
for $k=1,\ldots,n-1$. By Lemma~1, the processes $\XX_{(1)},\ldots,\XX_{(n-1)}$ have nondegenerate crossings as do $X_1,\ldots,X_n$, so the same holds for holds for $\XX_{(1)},\ldots,\XX_{(n-1)},X_n$.

Now, 
\[
X_{(1)}(t)=\XX_{(1)}(t)\lor X_n(t),\as,
\]
and we can apply Lemma~4, so by our inductive hypotheses,
{\allowdisplaybreaks
\begin{align*}
dX_{(1)}&=\1_{\{\XX_{(1)}(t)\ge X_n(t)\}}\dd\XX_{(1)}(t)+\1_{\{\XX_{(1)}(t)< X_n(t)\}}\dd X_n(t)\\
&=\1_{\{\XX_{(1)}(t)\ge X_n(t)\}}\sum_{i=1}^{n-1}\1_{\{X_i(t)=\XX_{(1)}(t)\}}\dd X_i(t)+\1_{\{X_n(t)=X_{(1)}(t)\}}\dd X_n(t)\\
&=\sum_{i=1}^{n-1}\1_{\{X_i(t)=X_{(1)}(t)\}}\dd X_i(t)+\1_{\{X_n(t)=X_{(1)}(t)\}}\dd X_n(t)\\
&=\sumi\1_{\{X_i(t)=X_{(1)}(t)\}}\dd X_i(t),\as
\end{align*}
}For $2\le k \le n-1$, we have
\[
X_{(k)}(t)= \XX_{(k-1)}(t)\land \big( \XX_{(k)}(t)\lor X_n(t) \big),\as
\] 
Since 
\[
\{t: \XX_{(k-1)}(t)= \big( \XX_{(k)}(t)\lor X_n(t) \big)\}\subset \{t: \XX_{(k-1)}(t)= \XX_{(k)}(t)\}\cup\{t: \XX_{(k-1)}(t)=X_n(t)\}
\]
and $d\brac{ \XX_{(k)}\lor X_n}_t\ll d\brac{ \XX_{(k)}}_t+d\brac{ X_n}_t$, it follows that  $\XX_{(k-1)}$ and 
$\XX_{(k)}\lor X_n$ have nondegenerate crossings. Hence, by Lemma~4, 
{\allowdisplaybreaks
\begin{align*}
dX_{(k)}(t)&=\1_{\{\XX_{(k-1)}(t)< \XX_{(k)}(t)\lor X_n(t)\}}\dd\XX_{(k-1)}(t)+
\1_{\{\XX_{(k-1)}(t)\ge \XX_{(k)}(t)\lor X_n(t)\}}\dd\big( \XX_{(k)}(t)\lor X_n(t)\big)\\
&=\1_{\{ \XX_{(k-1)}(t)< X_n(t)\}}\dd\XX_{(k-1)}(t)\\
&\qquad+\1_{\{\XX_{(k-1)}(t)\ge  X_n(t)\}}\big( \1_{\{\XX_{(k)}(t)\ge X_n(t)\}}\dd\XX_{(k)}(t)+\1_{\{\XX_{(k)}(t)< X_n(t)\}}\dd X_n(t)\big)\\
&=\1_{\{ \XX_{(k-1)}(t)<X_n(t)\}}\sum_{i=1}^{n-1}\1_{\{X_i(t)=\XX_{(k-1)}(t)\}}\dd X_i(t)\\
&\qquad+ \1_{\{\XX_{(k)}(t)\ge X_n(t)\}}\sum_{i=1}^{n-1}\1_{\{X_i(t)=\XX_{(k)}(t)\}}\dd X_i(t)+\1_{\{\XX_{(k-1)}(t)\ge X_n(t)>\XX_{(k)}(t)\}}\dd X_n(t) \\
&=\1_{\{ X_{(k)}(t)<X_n(t)\}}\sum_{i=1}^{n-1}\1_{\{X_i(t)=X_{(k)}(t)\}}\dd X_i(t)\\
&\qquad+ \1_{\{X_{(k)}(t)\ge X_n(t)\}}\sum_{i=1}^{n-1}\1_{\{X_i(t)=X_{(k)}(t)\}}\dd X_i(t)+\1_{\{X_n(t)=X_{(k)}(t)\}}\dd X_n(t) \\
&=\sum_{i=1}^{n-1}\1_{\{X_i(t)=X_{(k)}(t)\}}\dd X_i(t)+\1_{\{X_n(t)=X_{(k)}(t)\}}\dd X_n(t)\\
&=\sumi\1_{\{X_i(t)=X_{(k)}(t)\}}\dd X_i(t),\as
\end{align*}
}Finally, for $k=n$, we have
\begin{equation}
dX_{(n)}(t)=\sumi dX_i(t) -\sum_{k=1}^{n-1}dX_{(k)}(t),\as\tag*{\qed}
\end{equation}

\vspace{15pt}
\noindent{\large\bf Stratonovich representation for Atlas rank processes}
\vspace{10pt}
 
We would like to apply Proposition~1 to the Atlas model \eqref{0}. To do so we must show that the log-capitalization processes for an Atlas model are reversible and have nondegenerate crossings. 

\vspace{10pt}
\noindent{\bf Proposition~2.} {\em For the Atlas model \eqref{0}, the processes $\log X_1,\ldots,\log X_n$ are reversible and have nondegenerate crossings.}

\vspace{10pt}
\noindent{\em Proof.}  Girsanov's theorem and the properties of multidimensional Brownian motion imply that the processes $\log X_i$ of \eqref{0}  have nondegenerate crossings and that there are no {\em triple points,} i.e.,  for $i<j<k$, $\{t: \log X_i(t)=\log X_j(t)=\log X_k(t)\}=\empty$,~a.s.\ (see \citeN{Karatzas/Shreve:1991}). It remains to show that the $\log X_i$ are reversible.
 
Choose $k\in\{1,\ldots,n\}$ and  $t_0\in[0,T]$, and suppose that $\log X_j(t_0)=\log X_{(k)}(t_0)$. If for all $i\ne j$ we have $\log X_i(t_0)\ne \log X_j(t_0)$, then there is a neighborhood $U$ of $t_0$ in $[0,T]$ such that for $t\in U$, if $i\ne j$ then $\log X_i(t)\ne \log X_j(t)$. In this case, within $U$, the process $\log X_j$  is Brownian motion with drift, which is reversible.  Now suppose that $\log X_i(t_0)=\log X_j(t_0)$ for some $i\ne  j$. No-triple-points implies that there is a neighborhood $U$ of $t_0$ such that for $t\in U$, if $\ell\ne i,j$ then $\log X_i(t) \ne \log X_\ell(t)\ne \log X_j(t)$. Hence, within $U$ we can confine our attention to the two processes $\log X_i$ and $\log X_j$, in which case it was shown in \citeN{FIKP:2011} or \citeN{FIK:2013} that the time-reversed versions of these processes are continuous semimartingales. 

For each  $k$, the compactness of $[0,T]$ ensures that a finite subfamily of the neighborhoods $U$ will include all values of $t$, so the $\log X_i$  are reversible on  $[0,T]$. \qed

\vspace{10pt}
\noindent{\bf Corollary~1.} {\em For the Atlas model \eqref{0}, the rank processes $\log X_{(1)},\ldots,\log X_{(n)}$ satisfy
\begin{equation}\label{6.01}
d\log X_{(k)}(t)=\sumi\1_{\{X_i(t)=X_{(k)}(t)\}}\dd \log X_i(t),\as,
\end{equation}
}

\noindent{\em Proof.}  Follows immediately from Propositions~1 and~2. \qed

\vspace{15pt}
\noindent{\large\bf An application to portfolio return decomposition}
\vspace{10pt}
 
For $n\ge2$, consider a stock market of stocks with capitalizations represented by the positive continuous semimartingales $X_1.\ldots,X_n$ defined on $[0,T]$. The {\em market weight processes} $\m_1,\ldots,\m_n$ are defined by
\[
 \m_i(t)\eqdef\frac{X_i(t)}{X_1(t)+\cdots+X_n(t)},
\]
and the {\em ranked market weight processes} $\m_{(k)}$ are defined accordingly. If the processes $\log X_1,\ldots,\log X_n$ of a market are reversible and have nondegenerate crossings, the same will hold for the log-weight processes $\log\m_1,\ldots,\log\m_n$.

The {\em market portfolio} is the portfolio with weights $\m_i$ and portfolio value process
\[
Z_\m(t)=X_1(t)+\cdots+X_n(t),\as
\]
In \citeN{F:2016} it was shown that for a portfolio $\p$, the relative log-return  $d\log(Z_\p/Z_\m)$  can be decomposed into 
\[
d\log\big(Z_\p(t)/Z_\m(t)\big)=d\log {\mathcal S}_\p(t) + d\T_\p(t),\as,
\]
where ${\mathcal S}_\p$ is a {\em structural process}  defined by
\begin{equation}\label{7}
d\log{\mathcal S}_\p(t)\eqdef\sumi\p_i(t)\dd\log \m_i(t),
\vspace{-10pt}
\end{equation}
and  $\T_\p$ is a {\em trading process} with
\[
d\T_\p(t)\eqdef d\log\big(Z_\p(t)/Z_\m(t)\big)-d\log {\mathcal S}_\p(t),
\]
at least when the Stratonovich integrals in \eqref{7} are all defined.

Let $S$ be a real-valued $C^2$ function defined on a neighborhood of the unit simplex $\D^n\subset\R^n$. Then we shall say that the portfolio $\p$ is {\em generated by} the function $\S$ of the ranked market weights if $\S(\m(t))=S(\m_{(1)}(t),\ldots,\m_{(n)}(t))$ and the portfolio weight processes $\p_i$  are given by
\begin{equation}\label{7.1}
\p_{p_t(k)}(t)=\Big(D_k\log S(\m_{(\cdot)}(t))+1-\sumj\m_{(j)}(t)D_j\log S(\m_{(\cdot)}(t))\Big)\m_{(k)}(t),
\end{equation}
where  $p_t$ is the inverse of $r_t\in\Sigma_n$. In this case, the relative log-return of $\p$ will satisfy
 \[
d\log\big(Z_\p(t)/Z_\m(t)\big)=d\log\S(\m(t))+d\Theta(t),\as,
\]
where $\Theta$  is a function of locally bounded variation (see \citeN{F:rank} or \citeN{F:2002}, Theorem~4.2.1).

\vspace{10pt}
\noindent{\bf Proposition~3.} {\em Suppose that the market log-weight processes $\log\m_1,\ldots,\log\m_n$  are reversible and have nondegenerate crossings. Let $\p$ be the portfolio generated by the function $\S$ of the ranked market weights.  Then
\begin{equation}\label{100}
d\log{\mathcal S}_\p(t)=d\log\S(\m(t)),\as,
\end{equation}
and
\begin{equation}\label{101}
d\T_\p(t)=d\Theta(t),\as
\end{equation}
}

\vspace{-10pt}
\noindent{\em Proof.} By hypothesis, $\S(\m(t))=S(\m_{(\cdot)}(t))$, where $S$ is a real-valued $C^2$ function defined on a neighborhood of the unit simplex $\D^n\subset\R^n$. Then
{\allowdisplaybreaks
\begin{align}
d\log\S(\m(t))&=d\log S(\m_{(\cdot)}(t))\notag\\
&= \sumk D_k\log S(\m_{(\cdot)}(t))\dd\m_{(k)}(t)\notag\\
&=\sumk D_k\log S(\m_{(\cdot)}(t))\m_{(k)}(t)\dd\log\m_{(k)}(t)\notag\\
&=\sumk\p_{p_t(k)}(t)\dd\log\m_{(k)}(t)\label{7.5}\\
&=\sumk\p_{p_t(k)}(t)\sumi\1_{\{X_i(t)=X_{(k)}(t)\}}\dd\log\m_i(t)\label{8}\\
&=\sumi\sumk\p_{p_t(k)}(t)\1_{\{X_i(t)=X_{(k)}(t)\}}\dd\log\m_i(t)\notag\\
&=\sumi\p_i(t)\dd\log\m_i(t)\notag\\
&=d\log{\mathcal S}_\p(t)\as,\notag
\end{align}
}where \eqref{7.5} is due to \eqref{7.1} and the fact that
\[
\sumk\m_{(k)}(t)\dd\log\m_{(k)}(t)=\sumk d\m_{(k)}(t)=d\sumk\m_{(k)}(t)=0,\as,
\]
and \eqref{8} follows from Proposition~1. 

The representation for $\T_\p$ follows by construction. \qed

\vspace{10pt}
\noindent{\bf Corollary~2.} {\em For the Atlas model \eqref{0}, a portfolio generated by a function of the ranked market weights satisfies the decomposition \eqref{100} and \eqref{101}.}

\vspace{10pt}
\noindent{\em Proof.}  Follows immediately from Propositions~2 and~3. \qed

\vspace{10pt}
\noindent{\bf Remark.} The lemmata leading to Proposition~1 depend on the reversibility of  the semimartingales $X_i$. The localization argument in Proposition~2  to establish this reversibility for Atlas models depends on no-triple-points along with the $n=2$ results of \shortciteN{FIKP:2011} and \citeN{FIK:2013}. However, triple points may exist in {\em first-order models} and {\em hybrid Atlas models}, so  for these more general models localization to two dimensions fails and reversibility cannot immediately be established (see \shortciteN{BFK:2005}, \shortciteN{IPBKF:2011}, and \shortciteN{FIK:2012}). Hence, it appears that other methods may be needed to extend Corollary~2 to more general rank-based models.

\bibliographystyle{chicago}
\bibliography{math}

\end{document}